\documentclass[12pt]{article}

\usepackage{amsmath}
\usepackage{amssymb}
\usepackage{eucal}

\usepackage[english]{babel}

\begin{document}

\baselineskip 16pt

\title{ On the supersoluble residual\\ of mutually permutable products}


\author{V.\,S. Monakhov}


\maketitle

\begin{abstract}
We prove that if a group $G=AB$ is the mutually permutable product 
of the  supersoluble subgroups $A$ and~$B$, then
the supersoluble residual of $G$ coincides with 
the nilpotent residual of the derived subgroup $G^\prime$.
\end{abstract}



\bigskip

All groups in this paper are finite.

If $\mathfrak F$ is a formation and $G$ is a group, 
then $G^\mathfrak F$ is the $\mathfrak F$-residual of $G$, i.e.,
the smallest normal subgroup of $G$ with quotient in~$\mathfrak F$.
A group $G=AB$ is called the mutually permutable product
of subgroups $A$ and $B$ if $UB=BU$ and $AV=VA$ for all $U\le A$
and $V\le B$. Such groups were studied in~\cite{AS}--\cite{BH05},
see also~\cite{BEA}. 

We prove the following theorem.

\medskip

{\bf Theorem 1.} 
{\sl Let $G=AB$ be the mutually permutable product of 
the supersoluble subgroups $A$ and~$B$. 
Then $G^\mathfrak U=(G^\prime )^\mathfrak N=\big
{[}A,B\big {]}^\mathfrak N$.}

\medskip

Here $\mathfrak N$ and $\mathfrak U$ are respectively 
the formation of all nilpotent groups and
the formation of all supersoluble groups, and
$\big {[}A,B\big {]}=\langle [a,b]\mid  a\in A,b\in B\rangle$.

We need the following lemmas.

\medskip

{\bf Lemma 1} (\cite[4.8]{Mon}).
{\sl Let $G=AB$ be the product of two subgroups $A$ and~$B$.
Then

$(1)$ $\big {[}A,B\big {]}\lhd G$;

$(2)$ if $A_1\lhd A$, then $A_1\big {[}A,B\big {]}\lhd G;$

$(3)$ $G^\prime =A^\prime B^\prime \big {[}A,B\big {]}$.}

\medskip

If $\mathfrak X$ and $\mathfrak F$ are hereditary formations,
then, according to~\cite[p.~337-338]{DH}, the product
$$
\mathfrak X\mathfrak F=\{~G\in \mathfrak E \mid
G^{\mathfrak F}\in \mathfrak{X}\}
$$
is also a hereditary formation.

\medskip

{\bf Lemma 2} (\cite[IV.11.7]{DH}). {\sl
Let $\mathfrak F$ and $\mathfrak H$  be formations, 
$G$ be a group and $K\lhd G$. Then

$(1)$ $(G/K)^{\mathfrak F}=G^{\mathfrak F}K/K$;

$(2)$ $G^{\mathfrak F\mathfrak H}=(G^\mathfrak H)^\mathfrak F$;

$(3)$ if $\mathfrak H\subseteq \mathfrak F$, then
$G^\mathfrak F\subseteq G^\mathfrak H$.}

\medskip

If $H$ is a subgroup of a group $G$, then $H^G$ denotes 
the smallest normal subgroup of $G$ containing $H$.

\medskip

{\bf Lemma 3} (\cite[5.31]{Mon}). {\sl
Let $H$ be a subnormal subgroup of a group $G$.
If $H$ belongs to a Fitting class $\mathfrak F$,
then $H^G\in \mathfrak F$. In particular,

$(1)$ if $H$  is nilpotent, then $H^G$ is also nilpotent;

$(2)$ if $H$ is $p$-nilpotent, then $H^G$ is also $p$-nilpotent.}

\medskip

{\bf Lemma 4.}
{\sl
Let $G=AB$ be the product of the supersoluble subgroups $A$ and~$B$.
Then $G^\mathfrak U\le \big {[}A,B\big {]}$.}

\medskip

{\bf P r o o f.} 
By Lemma~1\,(1,3) and Lemma~2\,(1),
$$
\big {(}G/\big {[}A,B\big {]}\big {)}^\prime =G^\prime \big {[}A,B\big {]}/\big {[}A,B\big {]}=
A^\prime B^\prime \big {[}A,B\big {]}/\big {[}A,B\big {]}=
$$
$$
=\big {(}A^\prime \big {[}A,B\big {]}/\big
{[}A,B\big {]}\big {)}\big {(}B^\prime \big {[}A,B\big {]}/\big {[}A,B\big {]}\big {)}.
$$
The subgroups
$
\big {(}A^\prime \big {[}A,B\big {]}\big {)}/\big {[}A,B\big {]}\simeq
A^\prime /\big {(}A^\prime \cap \big {[}A,B\big {]}\big {)},
$
$$
\big {(}B^\prime \big {[}A,B\big {]}\big {)}/\big {[}A,B\big {]}\simeq 
B^\prime /\big {(}B^\prime \cap\big {[}A,B\big {]}\big {)}
$$
are nilpotent~\cite[VI.9.1]{Hup} and normal in~$G/\big{[}A,B\big {]}$ by Lemma~1\,(3),
so $\big {(}G/\big {[}A,B\big {]}\big {)}^\prime $ is nilpotent.
By Lemma~1\,(3), $A\big {[}A,B\big {]}$ and $B\big {[}A,B\big {]}$
are normal in~$G$. In view of the Baer Theorem~\cite{B57},
$G/\big {[}A,B\big {]}$ is supersoluble. Hence,
$G^\mathfrak U\le \big {[}A,B\big {]}$. \hfill $\boxtimes$

\medskip

A Fitting class which is also a formation is called a Fitting formation.
The class of all abelian groups is denoted by $\mathfrak A$.

\medskip

{\bf Lemma 5} \cite[II.2.12]{DH}.
{\sl Let $\mathfrak X$ be a Fitting formation, and let $G=AB$
be the product of normal subgroups $A$ and $B$.
Then $G^\mathfrak X=A^\mathfrak XB^\mathfrak X$.}

\medskip

{\bf P r o o f \ of \ Theorem~1.}
By Lemma~4, $G^\mathfrak U\le \big {[}A,B\big {]}$.
Since $\mathfrak U\subseteq \mathfrak N\mathfrak A$~\cite[VI.9.1]{Hup},
by Lemma~2\,(2,3), we have
$$
G^{(\mathfrak N\mathfrak A)}=(G^\mathfrak A)^\mathfrak N=
(G^\prime)^\mathfrak N\le G^\mathfrak U.
$$
Verify the reverse inclusion. Since 
$$
(G/(G^\prime )^\mathfrak N)^\prime =
G^\prime (G^\prime )^\mathfrak N/(G^\prime )^\mathfrak N =
G^\prime/(G^\prime )^\mathfrak N
$$
is nilpotent,
$$
G/(G^\prime )^\mathfrak N=A(G^\prime )^\mathfrak N/(G^\prime )^\mathfrak
N\cdot
B(G^\prime )^\mathfrak N/(G^\prime )^\mathfrak N
$$
is supersoluble in view of~\cite[Theorem 3.8]{AS} and $G^\mathfrak U\le (G^\prime )^\mathfrak N$.
Thus, $G^\mathfrak U= (G^\prime )^\mathfrak N$.

By Lemma~1\,(3),
$
G^\prime =A^\prime B^\prime [A,B]=(A^\prime)^G(B^\prime)^G [A,B].
$
The subgroups $A^\prime$ and $B^\prime$ are subnormal in~$G$ 
by~\cite[Theorem 1]{BH05} and nilpotent, therefore 
$(A^\prime)^G(B^\prime)^G$ is normal in~$G$ and nilpotent 
by Lemma~3\,(1). In view of Lemma~5 with $\mathfrak X=\mathfrak N$, 
we get
$$
G^\mathfrak U =(G^\prime )^\mathfrak N =((A^\prime)^G(B^\prime)^G)^\mathfrak N
[A,B]^\mathfrak N =[A,B]^\mathfrak N. \eqno \boxtimes
$$

\medskip

{\bf Corollary 1.1.}
{\sl Let $G=AB$ be the mutually permutable product of 
the supersoluble subgroups $A$ and~$B$.
If ${[}A,B\big {]}$ is nilpotent, then $G$ is supersoluble.}
\hfill $\boxtimes$

\medskip

The class of all $p$-nilpotent groups coincides with 
the product $\mathfrak E_{p^\prime}\mathfrak N_p$,
where $\mathfrak N_p$ is the class of all $p$\nobreakdash-\hspace{0pt}groups
and $\mathfrak E_{p^\prime}$ is the class of all $p^\prime$\nobreakdash-\hspace{0pt}groups.
A group $G$ is $p$-supersoluble 
if all chief factors of $G$
having order divisible by the prime $p$ are exactly of order $p$.
The derived subgroup of a $p$-supersoluble group is
$p$-nilpotent~\cite[VI.9.1\,(a)]{Hup}.
The class of all $p$\nobreakdash-\hspace{0pt}supersoluble groups 
is denoted by $p\mathfrak U$.
It's clear that
$
\mathfrak E_{p^\prime}\mathfrak N_p\subseteq p\mathfrak U\subseteq
\mathfrak E_{p^\prime}\mathfrak N_p\mathfrak A.
$

\medskip

{\bf Theorem 2.} {\sl
Let $G=AB$ be the mutually permutable product of the $p$-supersoluble
subgroups $A$ and~$B$. Then
$G^{p\mathfrak U}=(G^\prime)^{\mathfrak E_{p^\prime}\mathfrak N_p}=
\big {[}A,B\big {]}^{\mathfrak E_{p^\prime}\mathfrak N_p}$.}

\medskip

{\bf P r o o f.} By Lemma~2,
$$
(G^\prime)^{\mathfrak E_{p^\prime}\mathfrak N_{p}}=
(G^\mathfrak A)^{\mathfrak E_{p^\prime}\mathfrak N_{p}}=
G^{\mathfrak E_{p^\prime}\mathfrak N_{p}\mathfrak A}
\le G^{p\mathfrak U}.
$$
Verify the reverse inclusion. The quotient group
$$
G/(G^\prime )^{\mathfrak E_{p^\prime}\mathfrak N_p}=
(A(G^\prime )^{\mathfrak E_{p^\prime}\mathfrak N_p}/(G^\prime )^{\mathfrak
E_{p^\prime}\mathfrak N_p})
(B(G^\prime )^{\mathfrak E_{p^\prime}\mathfrak N_p}/(G^\prime )^{\mathfrak
E_{p^\prime}\mathfrak N_p})
$$
is the mutually permutable product of the $p$-supersoluble subgroups
$A(G^\prime )^{\mathfrak E_{p^\prime}\mathfrak N_p}/
(G^\prime )^{\mathfrak E_{p^\prime}\mathfrak N_p}$
and~$B(G^\prime )^{\mathfrak E_{p^\prime}\mathfrak N_p}/
(G^\prime )^{\mathfrak E_{p^\prime}\mathfrak N_p}$.
The derived subgroup
$$
(G/(G^\prime )^{\mathfrak E_{p^\prime}\mathfrak N_p})^\prime =
G^\prime (G^\prime )^{\mathfrak E_{p^\prime}\mathfrak N_p}/
(G^\prime )^{\mathfrak E_{p^\prime}\mathfrak N_p} =
G^\prime/(G^\prime )^{\mathfrak E_{p^\prime}\mathfrak N_p}
$$
is $p$-nilpotent. By~\cite[Corollary 5]{BH05}, 
$G/(G^\prime )^{\mathfrak E_{p^\prime}\mathfrak N_p}$ is $p$-supersoluble,
consequently, $G^{p\mathfrak U}\le (G^\prime )^{\mathfrak E_{p^\prime}\mathfrak N_p}$.
Thus, $G^{p\mathfrak U}= (G^\prime )^{\mathfrak E_{p^\prime}\mathfrak N_p}$.

By Lemma~1\,(3),
$
G^\prime =A^\prime B^\prime [A,B]=(A^\prime)^G(B^\prime)^G [A,B].
$
The subgroups $A^\prime$ and $B^\prime$ are subnormal 
in group~$G$~\cite[Theorem 1]{BH05} and $p$-nilpotent~\cite[VI.9.1\,(a)]{Hup},
hence $(A^\prime)^G(B^\prime)^G$ normal in~$G$ and
$p$-nilpotent by Lemma~3\,(2).
In view of  Lemma~5 with $\mathfrak X=\mathfrak E_{p^\prime}\mathfrak N_p$, 
we get
$$
G^{p\mathfrak U}= (G^\prime )^{\mathfrak E_{p^\prime}\mathfrak N_p}
=((A^\prime)^G(B^\prime)^G)^
{\mathfrak E_{p^\prime}\mathfrak N_p})[A,B]^{\mathfrak E_{p^\prime}\mathfrak N_p}=
[A,B]^{\mathfrak E_{p^\prime}\mathfrak N_p}. \eqno \boxtimes
$$

\medskip

{\bf Corollary 2.1.}
{\sl Let $G=AB$ be the mutually permutable product of 
the $p$-supersoluble subgroups $A$ and~$B$.
If ${[}A,B\big {]}$ is $p$-nilpotent, then $G$ is $p$-supersoluble.}

\bigskip


\noindent V.\,S. MONAKHOV

\noindent Department of mathematics, Gomel F. Scorina State
University, Gomel,  BELARUS

\noindent E-mail address: Victor.Monakhov@gmail.com


\bigskip

\begin{thebibliography}{19}

\bibitem{AS} 
M. Asaad and A. Shaalan, On the supersolubility of finite groups.
Arch. Math. 53, 318--326 (1989).

\bibitem{ABC}
M.\,J. Alejandre, A. Ballester-Bolinches and J. Cossey,
Permutable products of supersoluble groups.
J. Algebra 276, 453--461 (2004).

\bibitem{BCP}
A. Ballester-Bolinches, J. Cossey and M.\,C. Pedraza-Aguilera,
On products of supersoluble groups.
Rev. Mat. Iberoamericana 20, 413--425 (2004).

\bibitem{BH05}
J.\,C. Beidleman and H. Heineken,
Mutually permutable subgroups and group classes.
Arch. Math. 85, 18--30 (2005).

\bibitem{BEA}
A. Ballester--Bolinches, R. Esteban--Romero and M. Asaad,
Berlin--New York: Walter de Gruyter, 2010.
(De Gruyter expositions in mathematics;~53)

\bibitem{Mon}
V.\,S. Monakhov. Introduction to the theory of finite groups and their
classes. Minsk: Vyshejshaja shkola, 2006 (in Russian).

\bibitem{DH}
K. Doerk and T. Hawkes, Finite soluble groups.
Berlin--New York: Walter de Gruyter, 1992.

\bibitem{Hup}
B. Huppert. Endliche Gruppen I. Berlin -- Heidelberg -- New York: Springer. 1967.

\bibitem{B57}
R. Baer. Classes of finite groups and their properties.
Illinois J. Math. 1, 115-187 (1957).


\end{thebibliography}
\end{document}